\documentclass[11pt]{article}

\pdfoutput=1

\usepackage{latexsym}
\usepackage{amssymb}
\usepackage{amsfonts}
\usepackage{graphicx}

\def\T{{ \mathrm{\scriptscriptstyle T} }}

\usepackage{color}

\newcommand{\balpha}{\mbox{\boldmath$\alpha$}}

\newcommand{\bmu}{\mbox{\boldmath$\mu$}}
\newcommand{\bzero}{{\mathbf 0}}

\newcommand{\bA}{{\mathbf A}}

\newcommand{\bI}{{\mathbf I}}
\newcommand{\bM}{{\mathbf M}}
\newcommand{\bR}{{\mathbf R}}
\newcommand{\bT}{{\mathbf T}}
\newcommand{\bU}{{\mathbf U}}
\newcommand{\bV}{{\mathbf V}}
\newcommand{\bX}{{\mathbf X}}

\newcommand{\bZ}{{\mathbf Z}}

\newcommand{\bc}{{\mathbf c}}

\newcommand{\bm}{{\mathbf m}}
\newcommand{\bt}{{\mathbf t}}
\newcommand{\bu}{{\mathbf u}}
\newcommand{\bv}{{\mathbf v}}
\newcommand{\bx}{{\mathbf x}}
\newcommand{\by}{{\mathbf y}}
\newcommand{\bz}{{\mathbf z}}

\newcommand{\bSigma}{\mbox{\boldmath$\Sigma$}}

\newcommand{\bP}{{\mathbf P}}
\newcommand{\bS}{{\mathbf S}}
\newcommand{\bW}{{\mathbf W}}

\begin{document}

\vskip 1 cm

\vskip 1 cm

\smallskip

\begin{center}
{\LARGE
Measures of goodness of fit obtained by canonical transformations on Riemannian manifolds
}
\end{center}

\medskip

\begin{center}
 P.E. Jupp,  \\
{\small School of Mathematics and Statistics, University of St Andrews,  \\
St Andrews, Fife KY16 9SS, UK,}  \\
 \smallskip
   A. Kume, \\
{\small School of Mathematics, Statistics and Actuarial Science, University of Kent, Canterbury, Kent CT2 7NF, UK}
\end{center}
    
\medskip

\vskip 1cm

\begin{abstract}

The standard method of transforming a continuous distribution on the line to the
uniform distribution on $[0, 1]$ is the probability integral transform. 
Analogous transforms exist on compact Riemannian manifolds, $\mathcal{X}$, in that 
for each distribution with continuous positive density on  $\mathcal{X}$,
there is a continuous mapping of $\mathcal{X}$ to itself that transforms
the distribution into the uniform distribution. In general, this mapping is far from unique.
This paper introduces the construction of a version of such a probability integral that (under mild conditions) is canonical. 
The construction is extended to shape spaces,  simply-connected
spaces of non-positive curvature, and simplices.

The probability integral transform is used to derive tests of goodness of fit from tests of uniformity. 
Illustrative examples of these tests of goodness of fit are given 
involving (i) Fisher distributions on $S^2$, (ii) isotropic Mardia--Dryden distributions on the shape space $\Sigma^5_2$. 
Their behaviour is investigated by simulation.

\end{abstract}

\noindent
{\bf Keywords:}
Cartan--Hadamard manifold, Compositional data, Directional statistics, Exponential map, 
Probability integral transform, Shape space, Simplex.

\section{Introduction}

Directional statistics, shape analysis and compositional data analysis are concerned with probability distributions on Riemannnian manifolds, shape spaces and simplices, respectively.
The aim of this paper is to introduce and explore {a} canonical method of constructing transformations from such \linebreak
manifolds, 
$\mathcal{X}$, to certain associated manifolds, $\mathcal{Y}$, that send  \linebreak
arbitrary continuous distributions on $\mathcal{X}$ into standard distributions on $\mathcal{Y}$. 
More precisely, $\mathcal{Y}$ is $\mathcal{X}$ itself, or a tangent space to $\mathcal{X}$, or a star-shaped open  \linebreak
subset of a tangent space.
Given a basepoint $x$ in $\mathcal{X}$ and a standard  \linebreak
continuous distribution, $\nu$, on $\mathcal{Y}$, for any continuous distribution, $\mu$, on $\mathcal{X}$, we construct a function $\phi : \mathcal{X} \to \mathcal{Y}$ 
that is an almost-diffeomorphism (a diffeomorphism on the complement of some null set) that sends $\mu$ to $\nu$. 
Under mild conditions on uniqueness of medians of $\mu$ and of some distributions drived from it, 
$\phi$ as constructed here is canonical (in that any two versions differ only on a null set). 
These almost-diffeomorphisms, $\phi$, are used to obtain tests of goodness of fit to $\mu$ from tests of goodness of fit to $\nu$. 
If $\mathcal{X}$ is a compact Riemannian manifold then we can take 
$\mathcal{Y} = \mathcal{X}$, $\nu$ as the uniform distribution, and $\phi$ can be regarded as a form of probability integral transformation.
On compact manifolds our tests of goodness of fit complement the general 
Wald-type tests of Beran \cite{Beran79}, the score tests of Boulerice
and Ducharme \cite{BoulDuch97} and the Sobolev tests of Jupp \cite{Jupp05},
as well as the more specific tests in \cite{Kent82}, \cite{Mardiaetal84}, 
\cite{FisherBest84}, 
\cite{BestFisher86} (see \cite[Section 12.3]{MardiaJupp00}) 
and in \cite[Section 4.2]{Jonesetal15}, \cite{PewseyKato16} and 
\cite[Section 4.4]{Jupp15} in the case of copulae.

One important class of models on compact Riemannian manifolds $\mathcal{X}$ for which there is a canonical function
 $\phi : \mathcal{X} \to \mathcal{X}$ that takes any given distribution to the uniform distribution consists of those transformation models that are 
 obtained from the uniform distribution by the action of some group $G$ on $\mathcal{X}$, i.e., 
 the distribution of $x$ in $\mathcal{X}$ under parameter $g$ in $G$ is that of $gx$, where $x$ is uniformly distributed and $gx$ denotes 
 the image of $x$   \linebreak
 under $g$. Then $\phi (x)  = g^{-1} x$.  An example is the angular central Gaussian distributions on 
 the real projective space $\mathbb{R}P^{p-1}$ with probability density functions
\begin{equation}
f( \pm \bx;\bA) = | \bA | ^{- 1/2} (\bx ^{\top} \bA^{-1}\bx)^{-p/2} ,  \bx \in \mathbb{R}^p , 
\label{ACG}
\end{equation}
where $\bA$ is a non-singular symmetric $p \times p$ matrix and we may suppose that $| \bA | = 1$; see 
\cite[Section 9.4.4]{MardiaJupp00}.
Then $\phi (\pm \bx) =  \pm  \| \bA^{-1/2} \bx \|^{-1} \bA^{-1/2} \bx$.  
For general $\bA$, this transformation is different from that constructed in Section 2.2. If $\bA$ has only 2 distinct
eigenvalues then the two transformations are the same; see (b) after Remark 1.

The canonical transformations, $\phi$, are introduced in Section 2, first for spheres and then for compact Riemannian manifolds, 
shape spaces, Cartan--Hadamard manifolds and simplices. 
Section 3 shows how these transform\-ations send general tests of uniformity (or of goodness of fit to some standard distribution) into
general tests of goodness of fit.
The behaviour of these goodness-of-fit tests is illustrated in Section 4 by some simulation studies on the sphere, $S^2$, 
and on the shape space, $\Sigma^5_2$.

\section{Canonical transformations}

\subsection{Spheres}

Let $X$ be a random variable on the unit circle and suppose that an
orient\-ation and an initial direction on the circle have been chosen. Then the
{\it probability integral transformation} of the distribution is the 
transformation of the circle which sends $\theta$ to $U$, where
$U = 2 \pi \Pr (0 < X \le \theta)$.
If the distribution of $X$ is continuous then $U$ 
is distributed uniformly on the circle.
Thus  the probability integral transformation can be used to transform any test of   \linebreak
uniformity into a corresponding test of goodness of fit (see \cite[Section 6.4]{MardiaJupp00}).
For continuous distributions (with positive density) $\mu$, on $S^{p-1}$, the unit sphere in $\mathbb{R}^p$, with $p > 2$,
there are analogues 
$\phi : S^{p-1} \to S^{p-1}$ of the probability integral transformation 
that transform $\mu$ into the uniform distribution, $\nu$.
Such $\phi$ are far from unique, since if $\psi : S^{p-1} \to S^{p-1}$ preserves $\nu$ then the 
composite function $\psi  \circ \phi: S^{p-1} \to S^{p-1}$ also transforms $\mu$ into $\nu$.
Homeomorphisms $\psi$ that preserve $\nu$ can be constructed from embeddings 
$\gamma:  D^{p-1} \to S^{p-1}$ that map the uniform distribution on the disc, $D^{p-1}$,  to the uniform distribution on $\gamma (D^{p-1})$, 
together with functions $t \mapsto \bU_t$ from $[0,1]$ to the rotation group $SO(p-1)$ with $\bU_t = \bI_3$ for $t$ near $0$ or $1$. 
Then $\psi$ is the identity outside $\gamma (D^{p-1})$ and is given by $\psi \{ \gamma (r, \mathbf{\theta}) \} =  \gamma (r, \bU_r (\mathbf{\theta}))$ on 
$\gamma (D^{p-1})$, where $(r, \mathbf{\theta})$ are polar coordinates on $\gamma (D^{p-1})$.

Our construction of canonical versions of the probability integral transformation $\phi$ on $S^{p-1}$ is based on a set 
$S^{p-1} \supset S^{p-2} \dots  \supset S^s$ of nested spheres for which  
 \begin{equation}
\mbox{$S^{k-1}$ is the great sphere in $S^k$ normal to $\bm_k$ in $S^k$, for $k = p-1, \dots , s+1$,}  \\
\label{S^k}
 \end{equation}
where $\bm_k$ is some point in $S^k$.
The tangent-normal decomposition \cite[(9.1.20)]{MardiaJupp00}
 expresses each $\bx$ in $S^k$ as
\begin{eqnarray}
\bx &=& t \, \bm_k + \sqrt{1 - t^2} \, \bu 
\label{tan-nor1}  \\
&=& \cos (r) \, \bm_k + \sin (r) \, \bu ,  \label{tan-nor2}
\end{eqnarray}
where $t = \bx^{\top} \bm$, $\bu \in S^{k-1}$, the sphere normal to $\bm_k$, 
and $r = \arccos t$ is the colatitude of $\bx$. The function $p_k:  \bx  \mapsto \bu$ sends $S^k \setminus \{ \pm \bm_k \}$ into
$S^{k-1}$, so that, given a distribution $\mu$ on $S^{p-1}$, we can define distributions $\mu_{p-1}, \dots \mu_s$ on 
$S^{p-1}, S^{p-2} ,\dots  , S^s$ recursively by $\mu_{p-1} = \mu$ and
$\mu_{k-1}$ as the marginal distribution of $\bu$ on $S^{k-1}$ for $k = p-1, \dots , s + 1$.
We shall assume that
\begin{eqnarray}
 & & \mbox{$\mu$ is either uniform or has a unique (Fr\'echet) median $\bm_{p-1}$} , \qquad
\label{tunif.umimed1}  \\
 & & \mbox{for $k = p - 2, \dots s +1$, $\mu_k$ has a unique median $\bm_k$} ,
\label{tunif.umimed2}  \\
 & & \mbox{$\mu_s$ is the uniform distribution on $S^s$.} 
\label{tunif.umimed3}
\end{eqnarray}
If $\mu_1$ has a unique median $\bm_1$ then $\mu_0$ is automatically the uniform distribut\-ion on $S^0$.
The nested spheres in (\ref{S^k}) are reminiscent of the principal nested spheres of \cite{Jung.etal12}
but, whereas principal nested spheres may be small spheres and are chosen to give closest fit to the data, the spheres in
(\ref{S^k}) are great spheres and  are chosen to be orthogonal to $\bm_{p-1}, \dots , \bm_{s+1}$.
In cases in which (\ref{tunif.umimed1})--(\ref{tunif.umimed3}) hold, 
Proposition 1 provides a canonical version of the   \linebreak
probability integral transformation on $S^{p-1}$.

\bigskip

\noindent
{\bf  Proposition 1}

\emph{
Let $\mu$ be a probability distribution on $S^{p-1}$ such that the density of 
$\mu$ with respect to the  
uniform distribution, $\nu$, is continuous and positive.
Suppose that $\mu$ satisfies conditions (\ref{tunif.umimed1})--(\ref{tunif.umimed3}).
Then homeomorphic almost-diffeomorphisms $\phi_k : S^k \to S^k$ for $k = s , \dots , p - 1$ can be defined inductively by 
(a) $\phi_s$ is the identity,
(b) for $k = s + 1, \dots , p - 1$,
\begin{equation}
\phi_k (r, \bu) =   \psi_{k | \phi_{k-1} (\bu ) } (r)   \, \,\phi_{k-1} (\bu ) ,
\label{phi_k}
\end{equation}
where 
\[
\psi_{k | \bu } =    {\tilde F}_{\bu }^{-1} \circ F_{\bu} 
\]
with 
\begin{eqnarray}
 F_{\bu}(v) = & \Pr \left(  0 < R \le v | \bU = \bu \right) 
					\quad \mbox{under $\mu_k$}      \label{F_u}   \\
  {\tilde F}_{\bu}(v) = & \Pr \left(  0 < R \le v | \bU = \bu \right) 
					\quad \mbox{under $\nu_k$}    \label{Ftilde_u}
\end{eqnarray}
for $0 \le v \le \pi$, points $\bx$ in $S^{k+1}$ are identified with their coordinates $(r, \bu )$ as in (\ref{tan-nor2}), 
$(R, \bU)$ denotes a random element of $S^{k+1}$, and $\nu_k$ is the uniform distribution on $S^k$.
Then $\phi_{p-1}$ is a homeomorphic almost-diffeo\-morphism that} transforms $\mu$ into $\nu$.

\bigskip

\noindent
{\bf Proof}

From (\ref{phi_k}) and continuity of the density,  $\phi_k$ is a homeomorphism of $S^k$ and its restriction to 
$S^k \setminus \{ \pm \bm_k \}$ is a diffeo\-morphism.
It is straight\-forward to show that $\phi_{p-1}$ transforms $\mu$ into $\nu$.    \hspace{6 cm} $\Box$

\subsection{Compact Riemannian manifolds}

We now show how the probability integral transformation can be extended to arbitrary 
compact Riemannian manifolds in a canonical way.

Let $\mathcal{X}$ be a compact Riemannian manifold.
The Riemannian metric \linebreak
determines the volumes of infinitesimal cubes, and so equips $\mathcal{X}$ with a unique uniform probability measure, $\nu _{\mathcal{X}}$.
Let $\mu$ be a probability distribution on $\mathcal{X}$ having continuous positive density with respect to $\nu_X$.
If $\mathcal{X}$ is connected then there are homeomorphisms of $\mathcal{X}$ that 
transform $\mu$ into $\nu_X$; see \cite[Proposition 1]{Jupp15}.
One way of constructing such homeo\-morphisms, $\phi$, is by using the multivariate probability integral transformation 
(alias \emph{Rosenblatt transformation}, \cite{Rosenblatt52}) in coordinate neighbourhoods, as in the first proof in \cite{Moser65}. 
In the case in which the density is smooth, there is also a slick differential-geometric proof
 \cite[Theorem 2]{Moser65}.
This proof can be used to  provide a canon\-ical choice of $\phi$ but this involves solving a differential   
equation and does not give $\phi$ explicitly.
If $ \mathcal{X} = S^1$ or $\mathrm{dim} \, \mathcal{X} >1$ then, as in the spherical case, the homeo\-morphism $\phi$ is far 
from unique and it is not obvious how to make a canonical choice of  $\phi$. 
To obtain a canonical choice of $\phi$ by extending the construction in Proposition 1 to compact   \linebreak
Riemannian manifolds, 
we exploit the fact that, 
if $\mathcal{X} $ is a Riemannian manifold and $m$ is any point in $\mathcal{X}$ then the exponential map 
(see e.g., \cite[Section 1.6]{Helgason78})
from the tangent space, $T\mathcal{X} _m$,  at $m$ into $\mathcal{X}$ defines
 a system of Riemann\-ian \emph{normal coordinates} around $m$ as follows.
The inverse of this coordinate system maps the open set
$\left\{ (r, \bu) : 0 \le r < r_{\bu} , \bu \in T_1\mathcal{X} _m \right\}$ 
diffeomorphically onto an open set $\mathcal{B}$ of $\mathcal{X}$ by 
\begin{equation}
(r, \bu ) \mapsto \exp (r \bu) ,
\label{normcoords}  
\end{equation} 
where $T_1\mathcal{X} _m$ denotes the set of unit tangent vectors at $m$ and 
\[
r_{\bu} = \sup \{ r : \mbox{there is a unique minimising geodesic from $m$ to $\exp (r \bu)$} \} .
\] 
For $\mathcal{X} = S^{p-1}$, 
the tangent-normal decomposition (\ref{tan-nor1}) is related to the  \linebreak
normal coordinates by 
$t = \cos r$.
If $\mathcal{X}$ is compact then $\mathcal{X} \setminus \mathcal{B}$ has measure zero.
See, e.g.,    \cite[Proposition 2.113, Corollary 3.77, Lemma 3.96]{Gallotetal93}.
Thus absolute\-ly continuous probability distributions on $\mathcal{X}$ can be identified with absolutely continuous probability distributions 
on $\left\{ (r, \bu) : 0 \le r < r_{\bu} , \bu \in T_1\mathcal{X} _m \right\}$.
In particular, such a distribution induces a marginal distribution on 
$T_1\mathcal{X} _m$. 

\bigskip

\noindent
{\bf  Proposition 2}

\emph{
Let $\mu$ be a probability distribution on a compact Riemannian   \linebreak
manifold $\mathcal{X}$ of dimension $d$ such that the density of 
$\mu$ with respect to the uniform distribution, $\nu$, is continuous and positive.
Suppose that $\mu$ is   \linebreak
either uniform or has a unique median, $m$.  
If $\mu$ is the uniform distribut\-ion then define     
$\phi: \mathcal{X} \to \mathcal{X}$ as the identity.   
If $\mu$ is non-uniform then let  
 $\left\{ (r, \bu) : 0 \le r < r_{\bu} , \bu \in T_1\mathcal{X} _m \right\}$ be (maximal) Riemannian normal co\-ordinates on $\mathcal{B}$ 
with $m$ corresponding to the origin.
Assume that the marginal distributions on $T_1 \mathcal{X}_m$ obtained from $\mu$ and $\nu$ by using (\ref{normcoords}) satisfy conditions 
(\ref{tunif.umimed1})--(\ref{tunif.umimed3}).}

\emph{
Define the function $\phi : \mathcal{X} \to \mathcal{X}$ by
\begin{equation}
\phi \left\{ \exp ( r\bu ) \right\} 
=   \exp  \left[ 
{\tilde F}_{\psi _{d-1} ( \bu)}^{-1} \left\{ F_{\bu}(r) \right\} \,  \psi_{d-1}( \bu )
                   \right] \qquad  r \bu \in 
     \exp ^{-1}( \mathcal{B} )
\label{phi}
\end{equation}
and arbitrarily on $\mathcal{X} \setminus \mathcal{B}$, where 
$F_{\bu}$ and ${\tilde F}_{\bu}$ are defined by (\ref{F_u}) and (\ref{Ftilde_u}),
$\psi_{d-1} =  \tilde{\phi}_{d-1}^{-1} \circ \phi_{d-1}$ with 
$\phi_{d-1}, \tilde{\phi}_{d-1} : T_1\mathcal{X} _m \to T_1\mathcal{X} _m$ being the canonical uniformising 
almost-diffeo\-morphisms corresponding to $\mu$ and $\nu$, respectively,
given by Proposition 1 and identification of $T_1\mathcal{X} _m$ with $S^{d-1}$.
Then $\phi$ is a diffeo\-morphism almost everywhere and transforms $\mu$ into $\nu$.
}      

\medskip
\noindent
{\bf  Proof}

This is a straight\-forward calculation.   \hspace{5 cm} $\Box$

\medskip
\noindent
{\bf  Example}

\emph{The torus, $S^1 \times S^1$, can be written as $[- \pi, \pi] \times [- \pi, \pi]$, where $- \pi$ and $\pi$ are 
identified. Then $\mathcal{B}$ can be taken as
$(- \pi , \pi)  \times  (- \pi , \pi)$ and 
\[
  r_u = \pi / \max  \{  | \cos(u) |, |  \sin(u) |  \}   \qquad u \in [0, 2 \pi] . 
\] }

\medskip

We call the almost-diffeomorphism $\phi$ of Propositions 1 or 2 the 
\emph{probability integral transformation}.
It is canonical, since it is determined 
(except on null sets) by unique medians at each stage.

\bigskip

\noindent
{\bf  Remark 1}

\emph{The appropriate general mathematical setting for the constructions in Propositions 1 and 2 is that of 
ortho\-normal frames in a tangent space.
An ortho\-normal frame at a point $m$ in a $d$-dimensional manifold $\mathcal{X}$ is an   \linebreak
ordered set of ortho\-normal vectors in the tangent space $T\mathcal{X}_m$.
Let $\mu$ be a probability distribution on $\mathcal{X}$ such that the density of $\mu$ with respect to the   
uniform distribution, $\nu$, is continuous and positive.
Let $(m_{d-1}, \ldots , m_{s+1})$ be an ortho\-normal frame at $m$ 
and suppose that the distribution on the $u$-sphere normal to $m_{d-1}, \ldots , m_{s+1}$ is uniform. 
Then replacing the success\-ive medians in Propositions 1 and 2 by $m, m_{d-1}, \ldots , m_{s+1}$ 
defines an almost-diffeomorphism $\phi$ of $\mathcal{X}$ that takes $\mu$ to $\nu$.}
\hspace{5 cm} $\Box$

\medskip

A class of distributions for which the probability integral transformation takes a particularly simple form 
consists of those 
with unique median $m$ on $\mathcal{X}$ and for which the corresponding marginal distribution on  $T_1 \mathcal{X}_m$
(obtained using (\ref{normcoords})) is uniform.
If $\mathcal{X}$ is the sphere $S^{p-1}$, 
the projective space $\mathbb{R}P^{p-1}$, 
the rotation group $SO(3)$ 
or the complex project\-ive space $\mathbb{C}P^{k-2}$, then
these (include the distributions that have rotational symmetry about the unique median. 
Some examples are:

\begin{enumerate}
\item[(a)]   
For a distribution $\mu$ on $S^{p-1}$ that is rotationally symmetric about a unit vector $\bmu$,
the transformation $\phi$ given by (\ref{phi}) 
that sends $\mu$ into the uniform distribution has the form
\begin{equation}
 \phi (\bx) =  u \bmu +  \sqrt{ (1 - u^2)/(1 - t^2)} \left( \bI_p - \bmu \bmu ^{\top} \right) \bx , 
\label{phi.sym}
\end{equation}
where $t = \bx^{\top} \bmu$,
$\bI_p$ denotes the $p \times p$ identity matrix and
\begin{equation}
    u = G_0^{-1} (G_{\mu} (t)) ,
\label{Ginv-G}
\end{equation}
$G_{\mu}$ and  $G_0$ denoting the cumulative distribution functions of $\bx^{\top} \bmu$ 
when $\bx$ has  distribution $\mu$ and the uniform distribution, respectively.
In particular, for the Fisher distribution, $F(\bmu , \kappa)$, on $S^2$ 
with mean direct\-ion $\bmu$ and concentration $\kappa$, 
\begin{equation}
 u = \left(  2 e^{\kappa t}  - e^{\kappa} - e^{- \kappa}  \right) /
\left( e^{\kappa} - e^{- \kappa} \right) , \qquad  \kappa > 0 
\label{u.Fisher}
\end{equation} 
and $u = t$ for $\kappa = 0$ (see \cite[Example 1]{Jupp15}).

\item[(b)] Among the angular central Gaussian distributions on 
the real project\-ive space $\mathbb{R}P^{p-1}$ 
with probability density functions (\ref{ACG}),
those that are symmetrical about the modal axis $\pm \bmu$ have 
$\bA = a \, \bmu \bmu^\top+ b \, ( \bI_p -   \bmu \bmu^\top)$ with $a > b > 0$.
Then $\phi$  is given by     
\[
 \phi (\pm \bx) =  
\pm \left\{  u \bmu +  \sqrt{ (1 - u^2)/(1 - t^2)} \left( \bI_p - \bmu \bmu ^{\top} \right) \bx  \right\}, 
\]
where $t = \bx ^\top \bmu$ and 
\[
u = t / \sqrt{a/b + (1 - a/b) t^2} .
\] 
The transformation $\phi$ coincides with the standard transformation   \linebreak
$\pm \bx \mapsto \pm \| \bA^{-1/2}  \bx\|^{-1} \bA^{-1/2} \bx$ to uniformity on $\mathbb{R}P^{p-1}$
\cite[Section 9.4.4]{MardiaJupp00}.

\item[(c)] For the matrix Fisher distribution on $SO(3)$
with density proportional to $\exp  \left\{ \mathrm{tr} \left(\kappa  \bX^{\top} \bM \right) \right\}$ 
for $\kappa \ge 0$ and $\bM$ in $SO(3)$,
\cite[Example 2]{Jupp15} shows that 
$\bM^{\top} \bX$ and $\bM^{\top} \phi (\bX)$ have the same rotation axis,
and that the rotation angle, $u$, of $\bM^{\top} \phi (\bX)$ is related to  the rotation angle, $t$, of $\bM^{\top} \bX$ by 
\[
 {\tilde F}_0 (u)/ {\tilde F}_0 (\pi) =  {\tilde F}_{\kappa}(t)/ {\tilde F}_{\kappa}(\pi)  ,
\]
where 
${\tilde F}_{\kappa}(\theta) 
=  \int_0^{\theta}  e^{4 \kappa \cos ^2 (\omega /2)} \sin ^2 (\omega /2) d\omega$.
\item[(d)]  
On the shape space 
$\Sigma ^k_2$ of $k$ non-identical labelled landmarks in $\mathbb{R}^2$,  
the isotropic Mardia--Dryden distributions, alias isotropic offset normal distributions,
$MD( [ \bmu ], \kappa)$ \cite[Section 11.1.2]{DrydenMardia16}
of shapes $[\bX]$ obtained by isotropic Gaussian perturbation of the landmarks of shapes $[\bmu]$ have densities 
\begin{equation}
 f([\bX] ;  [ \bmu ], \kappa ) = 
 e^{\kappa \{ 1-\cos^2\rho([\bX], [\bmu]) ]\} } 
 \mathcal{L}_{k-2}( -  \kappa \cos^2\rho([\bX], [\bmu])) ,
\label{KVMILD}
\end{equation}
 where $\mathcal{L}_{k-2}$ is the Laguerre polynomial of order $k-2$, $\rho$ is the Riemann\-ian shape 
distance and $\kappa$ is a concentration parameter \cite[equations (11.11),(11.15)]{DrydenMardia16}.
Identification of $2 \times (k-1)$ real matrices $\bZ$ satisfying 
$\mathrm{trace} (\bZ \bZ  ^{\top}) =1$ with unit vectors $\bz$ in $\mathbb{C}^{k-1}$
leads to identification of the  space $\Sigma ^k_2$ with the complex projective space $\mathbb{C}P^{k-2}$.
Calculation shows that for the distribution with density (\ref{KVMILD}), the homeo\-morphism $\phi$ is 
\[
\phi ([\bz]) = [ u  \bmu + \sqrt{(1 - u^2)/(1 - t^2}) \left\{ \bz - (\bz^{\top} \bmu ) \bmu \right\}  ],
\]
where $t = \cos \, \rho( [X], [\mu])$, $u^2 = F_{[X], 0}^{-1} \left\{ F_{[X],\kappa} (t^2) \right\}$ with
$F_{[X], \kappa}$ defined by
\begin{eqnarray*}
  & & F_{[X], \kappa} ( x )   \nonumber  \\  
    &=&  (k-2) e^{ \kappa}  \sum_{i=0}^{k-2} \sum_{r=0}^{k-3} {k-2 \choose i} {k-3 \choose r}
\frac{ (-1)^r \kappa^i} {i!}   
\int_{0}^x e^{ - \kappa s} s^{r+i} d s   .    \qquad 
\label{F_Xkappa}
 \end{eqnarray*}  
For $\kappa = 0$ (corresponding to the uniform distribution)
$F_{[X], \kappa}$ takes the simple form
\[
  F_{[X],0}(x)  = 1-(1-x)^{k-2}  . 
 \]
\end{enumerate}

\subsection{Shape spaces}

The probability integral transformation can be defined also for the shape spaces, $\Sigma_m^k$, of shapes of $k$ non-identical 
labelled landmarks in $\mathbb{R}^m$. 
As indicated after (\ref{KVMILD}), the
 space $\Sigma_2^k$ can be identified with the complex projective space $\mathbb{C}P^{k-2}$, and so is a compact Riemannian manifold. 
For $m > 2$, $\Sigma_m^k$ is not a manifold but for our purposes, it is enough to work on the \emph{non-singular part} 
of $\Sigma_m^k$, which is the open set consisting of the shapes of $k$ non-identical 
labelled landmarks in $\mathbb{R}^m$ that do not lie in any $(m-2)$-dimensional affine subspace.

It follows from \cite[Section 6.3 and Theorem 6.5]{Kendalletal99} that, 
for $x$ in  the non-singular part of $\Sigma_m^k$ there is a system of Riemannian normal coordinates with inverse that maps 
an open set $\left\{ (r, \bu) : 0 \le r < r_{\bu} , \bu \in T_1\mathcal{X} _x \right\}$ 
diffeomorphically onto an open set $\mathcal{B}$ of $\Sigma_m^k$ by (\ref{normcoords}),  
where $T_1\mathcal{X} _x$ denotes the set of unit tangent vectors at $x$, and 
$\Sigma_m^k \setminus \mathcal{B}$ has measure zero.
If the  distribution on $T_1\mathcal{X}$ 
satisfies conditions (\ref{tunif.umimed1})--(\ref{tunif.umimed3}) then the probability integral transform can be defined as in Proposition 2.

\subsection{Cartan--Hadamard manifolds}

The Cartan--Hadamard manifolds are the complete simply-connected mani\-folds with non-positive curvature.
It follows from the Cartan--Hadamard theorem 
\cite[Theorem I 13.3]{Helgason78}, \cite{KobayashiNomizu69}
that on a Cartan--Hadamard manifold, $\mathcal{X}$,  
the inverse of the exponential map at any basepoint $x$ identifies $\mathcal{X}$ with $T\mathcal{X}_x$.
Then the choice of some `basepoint' distribution $\nu$ on $\mathcal{X}$ enables an extension of the approach used in Section 3. 
Important instances of such manifolds are the simplicial shape spaces   
 of shapes of $(m+1)$-simplices in $\mathbb{R}^m$ with positive volume, 
 equipped with a Riemannian metric derived from a natural metric on $SL(m)$  \cite[Section 3.6.2]{Small96}, \cite[Section 3]{Le.Small99}.
The case $m=2$ gives the space of shapes of non-degenerate triangles in the plane, which can be identified with the
Poincar\'e half-plane, $\mathbb{H}^2 = \{ (x_1, x_2) \in \mathbb{R}^2 : x_2 > 0 \}$, with Riemannian metric
$g_{ij }= \delta_{ij} x_2^{-2}$.
This space was used in \cite{Huckemannetal10} as a sample space for electrical impedances.

\medskip
\noindent
{\bf  Proposition 3}

\emph{
Let $\mu$ and $\nu$ be probability distributions on a Cartan--Hadamard manifold, $\mathcal{X}$.
Let $m$ be a point of $\mathcal{X}$ and $\left\{ (r, \bu) : 0 \le r , \bu \in T_1\mathcal{X} _m \right\}$ be Riemann\-ian normal coordinates on 
$\mathcal{X}$ with $m$ corresponding to the origin.
 Define the function $\phi : \mathcal{X} \to \mathcal{X}$ by
\[
  \phi \left\{ \exp ( r\bu ) \right\} 
=   \exp  \left[ 
    {\tilde F}_{\psi _{d-1} ( \bu)^{-1} \left\{ F_{\bu}(r) \right\} \,  \psi_{d-1}( \bu ) } \right] ,
\]
where $F_{\bu}$ and ${\tilde F}_{\bu}$ are defined by (\ref{F_u}) and (\ref{Ftilde_u}),
and $\psi_{d-1} =  \tilde{\phi}_{d-1}^{-1} \circ \phi_{d-1}$ with 
$\phi_{d-1}, \tilde{\phi}_{d-1} : T_1\mathcal{X} _m \to T_1\mathcal{X} _m$ being the canonical uniformising 
almost-diffeo\-morphisms corresponding to $\mu$ and $\nu$, respectively.}

\emph{
Then $\phi$ is an almost-diffeo\-morphism that
maps geodesics through $m$ into geodesics through $m$ and transforms $\mu$ into $\nu$.}

\subsection{Simplices}

The open $(p-1)$-simplex is
\[
\Delta_{p-1} = \left\{ (y_1, \ldots , y_p) :     y_j > 0, \, \sum_{j=1}^p y_j = 1\right\} .
\]
There is a canonical base point, the centroid, $\bc = (p^{-1}, \ldots p^{-1})$ and a canonical Riemannian metric 
obtained by regarding $\Delta_{p-1}$ as an affine subspace of $\mathbb{R}^p$.
The unit tangent sphere at $\bc$ is 
\[
 T_1 \Delta_{p-1, \bc}
=  \left\{ \bu = (v_1, \ldots , v_p) : \sum_{j=1}^p v_j = 0,   \sum_{j=1}^p v_j ^2= 1  \right\}
\]
and the exponential map is
 \begin{equation}
\exp ( r \bu ) = \bc + r \bu 
\label{exp.simplex}
\end{equation}
for $r \in [0, 1/(p \max_{1 \le j \le p} | v_j |) ]$.
The uniform distribution is a scaled version of Lebesgue measure 
on $\Delta_{p-1}$
and the corresponding marginal distribution on the unit tangent sphere is the uniform distribution on $T_1 \Delta_{p-1, \bc}$.

\subsubsection{Using the exponential map}

The manifold $\Delta_{p-1}$ is simply connected and has curvature 0 but it is not complete.
The exponential map (\ref{exp.simplex}) is a diffeo\-morphism between a star-shaped portion of 
$T \Delta_{p-1, \bc}$ and $\Delta_{p-1}$.
Let $\mu$ be a distribution on $\Delta_{p-1}$ with continuous positive density
with respect to the uniform  distribution, $\nu$.
Then a minor variant of Proposition 3 produces a canonical almost-diffeomorphism
 $\phi : \Delta_{p-1} \to \Delta_{p-1}$ that transforms $\mu$ into $\nu$.

\medskip

\noindent
{\bf  Proposition 4}

\emph{
Let $\mu$ be a probability distribution on $\Delta_{p-1}$ having continuous positive density with respect to Lebesgue measure.
Let ${\bc}$  be the barycentre of $\Delta_{p-1}$
and $\left\{ (r, \bu) : 0 \le r , \bu \in T_1 \Delta_{p-1 , \bc}  \right\}$ be Riemannian normal coordinates on $\Delta_{p-1}$ with ${\bc}$ 
corresponding to the origin.
Let $(R, \bU)$ be the  normal coordinates of a random element of $\Delta_{p-1}$.
Define the function $\phi: \Delta_{p-1} \to \Delta_{p-1}$ by
\[
\phi \left\{ \exp ( r\bu ) \right\} 
= \exp  \left[  {\tilde F} _{\psi (\bu )}^{-1} \left\{ F_{\bu}(r) \right\} \,  \psi (\bu ) \right]  ,
\]
where 
$F_{\bu}$ and ${\tilde F}_{\bu}$ are defined by (\ref{F_u}) and (\ref{Ftilde_u})
and $\psi : T_1 \Delta_{p-1, \bc} \to T_1 \Delta_{p-1, \bc}$ is the almost-canonical homeo\-morphism such that $\psi (\bU)$ is uniformly distributed.
Then $\phi$ is a diffeo\-morphism almost everywhere,
maps geodesics through ${\bc}$ into geodesics through ${\bc}$, and transforms $\mu$ into $\nu$.}
\label{prop: prop4}

\subsubsection{Using radial projection}

An alternative to using the exponential map (\ref{exp.simplex}) is to use `radial  \linebreak projection' of $\Delta_{p-1} \setminus \{ \bc \}$ onto its boundary $\partial \Delta_{p-1}$.
The coordinates $(r, z_1, \ldots , z_p)$ given by radial projection are defined by
\begin{eqnarray}
 r &=& 
\left\{\begin{array}{lll}
           0       & \mbox{if} & \bx = \bc ,\\
 1 - p y_{(1)} & \mbox{if} & \bx \ne \bc 
\end{array} \right.  \label{rad.proj.r}   \\
  z_j &=& r^{-1} (y_j - y_{(1)}) \quad  j = 1, \dots , p ,   \label{rad.proj.z} 
\end{eqnarray}
$y_{(1)}$ denoting the smallest of $y_1, \ldots , y_p$.
Then $r \in [0, 1$).
A simple calculation shows that the density of the uniform distribution with respect to $dr \, dz_1 \, \ldots  \, dz_{i-1} \, dz_{i+1} \, \ldots \, dz_p$ is proportional to
$r^{-(p-1)} $.      It follows that, for $i = 1, \ldots, p$, 
radial projection of $\Delta_{p-1,i} = \left\{ (y_1, \ldots , y_p) \in \Delta_{p-1} \setminus \{ \bc \} :   y_{(1)} =   y_i \right\}$ 
onto the face $\partial _i \Delta_{p-1} = \left\{ (z_1, \ldots , z_p) :  z_i = 0 \right\}$ sends the uniform distribution on $\Delta_{p-1,i}$ to 
the uniform distribution on the $(p-2)$-simplex $\partial _i \Delta_{p-1}$.
The boundary, $\partial \Delta_{p-1}$, of $\Delta_{p-1}$ is the union of $\partial _1 \Delta_{p-1}, \ldots, \partial _p \Delta_{p-1}$.

The next proposition shows that radial projection provides canonical uniformising homeomorphic almost-diffeo\-morphisms of simplices that are analoguous  to 
those for spheres that are described in Proposition 1. \linebreak
Unlike the construction in Proposition 1, the  construction in Proposition 5 does not assume uniqueness of medians, as in (\ref{tunif.umimed1})--(\ref{tunif.umimed2}).

\medskip

\noindent
{\bf  Proposition 5}

\emph{
Let $\mu$ be a probability distribution on $\Delta_{p-1}$ having continuous positive density with respect to Lebesgue measure.
For $k = 0, \dots, p-2$, denote by  $\partial^{p-1-k}\Delta_{p-1}$, the union of the  $k$-dimensional faces of $\Delta_{p-1}$.
Then repeated radial projection sends $\mu$ to a probability distribution $\mu_k$ on $\partial^{p-1-k}\Delta_{p-1}$.
Let $s$ be the largest value of $k$ for which $\mu_k$ is uniform.
For $k = s + 1, \dots , p - 1$, let $r, z_1, \dots, z_{i-1}, z_{i+1}, \dots z_k$ be coordinates 
(defined analogously to those in (\ref{rad.proj.r})--(\ref{rad.proj.z})) on the part of
the $(p-1-k)$-simplex in $\partial^{p-1-k}\Delta_{p-1}$ on which  $z_i = 0$.
Define functions $\phi_k : \partial^{p-k-1}\Delta_{p-1} \to \partial^{p-k-1}\Delta_{p-1}$ for $k = s , \dots , p - 1$ recursively by 
(a) $\phi_s$ is the identity,
(b) for $k = s + 1, \dots , p - 1$,
\begin{equation}
\phi_k (r, \bz) =   F_{\bz}(r) ^{1/(k+2-p)}   \, \,\phi_{k-1} (\bz ) ,
\label{phi_k}
\end{equation}
where 
\begin{equation}
 F_{\bz}(v) =  \Pr \left(  0 < R \le v | \bZ = \bz \right) 
					\quad \mbox{under $\mu_k$}      \label{F_u}   
\end{equation}
for $0 \le v \le 1$ and $\bz = (z_1, \dots, z_{i-1}, z_{i+1}, \dots z_k)$,
Then $\phi_{p-1}$ is a homeomorphic almost-diffeo\-morphism that transforms $\mu$ into $\nu$.}
\label{prop: prop5}

\bigskip

\noindent
{\bf Proof}

This is a straightforward calculation using the fact that   \linebreak
$\Pr \left(  0 < R \le v | \bZ = \bz \right) = v^{k+2-p}$ under the uniform distribution on this   \linebreak
$(p-1-k)$-simplex.    \hspace{9 cm} $\Box$

\medskip
\noindent
{\bf  Remark 2}

\emph{
The almost-canonical homeomorphisms $\phi$  introduced in 
this Section can be used in the simulation of 
arbitrary continuous distributions on $\mathcal{X}$. 
Let $\mu$ and $\nu$ be probability distributions on $\mathcal{X}$ 
and $\phi$ any transformation that takes $\mu$ into $\nu$. 
If $x_1, \ldots , x_n$ in $\mathcal{X}$ are a random sample from  $\nu$ then  
$\phi ^{-1}(x_1), \ldots , \phi ^{-1} (x_n)$ are a random sample from  $\mu$.}

\section{Goodness-of-fit tests via transformation}

Let $\mu$ and $\nu$ be probability distributions on $\mathcal{X}$. 
Then any transformation, $\phi$, that takes $\mu$ into $\nu$ 
can be used to transform any test, $T$, of goodness of fit to $\nu$ into a test, $\phi ^* T$, of goodness of fit to $\mu$. 
Given points $x_1, \ldots , x_n$ in $\mathcal{X}$, $\phi ^* T$ is obtained by applying $T$ 
 to the transformed data, $\phi (x_1), \ldots , \phi (x_n)$.
The null distribution of {$\phi ^* T$ is the same as that of $T$.

Often the null hypothesis about the distribution generating the data is not that it is some specified distribution but that it is a distribution in a given 
parametric model, $\left\{ \mu_{\theta}  : \theta \in \Theta \right\}$.   
For each $\theta$ in $\Theta$, let $\phi_{\theta}$ be a transformation that takes $\mu_{\theta}$ into $\nu$. 
Let ${\hat \theta}$ be an estimate of $\theta$.
Then goodness of fit to $\left\{ \mu_{\theta}  : \theta \in \Theta \right\}$ is tested by applying 
$T$ to the  transformed data, $\phi _{{\hat \theta}} (x_1), \ldots , \phi _{{\hat \theta}} (x_n)$.
Significance can be assessed by simulation from the fitted distribution.
If a good approximation to the null distribution of $T$ is available then simulation can be avoided by using this approximation.

Provided that the estimator giving ${\hat \theta}$ is consistent, the consistency properties of 
$\phi ^* T$ 
are inherited from those of $T$.   
In particular, if ${\hat \theta}$ is the maximum likelihood estimate then $\phi ^* T$   
  is consistent against all alternatives if and only if $T$   
is consistent against all alternatives.
}

\subsection{Spheres}

On a sphere the uniform distribution provides a canonical choice for $\nu$. Then the transformation, $\phi$, of Proposition 1
 that takes $\mu$ into $\nu$ can be used to transform tests of uniformity into tests of goodness of fit to $\mu$.

One nice characterisation of the uniform distributions on $S^2$ is that, 
for a uniformly distributed random vector with longitude $\psi$ and colatitude $\theta$, 
(a) $\psi$ is uniformly distributed on $[0, 2 \pi ]$, 
(b) $\cos \theta$ is uniformly distributed on $[-1 , 1]$,
(c)  $\psi$ and $\theta$ are independent.
Thus combining any tests of (a), (b) and (c) gives a test of uniformity on $S^2$.
Using the general construction given in the previous paragraph with $\phi :S^2 \to S^2$ given by  (\ref{phi}) but with (\ref{u.Fisher}) 
replaced by the approximation $2 e ^{\kappa (t-1)} -1$ to (\ref{u.Fisher}) for $\kappa$ not close to $0$,
taking the tests in (a), (b) and (c) to be Kuiper's $V_n$, the Kolmogorov--Smirnov test, and a rather special `2-variable' test yields the 
standard method \cite[Section12.3.1]{MardiaJupp00} of investigating goodness of fit of  Fisher distributions on $S^2$.

\subsection{Compact Riemannian manifolds and shape spaces}

On a compact Riemannian manifold or a shape space the uniform distribution provides a canonical choice for $\nu$. Then the transformation, $\phi$, of Proposition 2 that takes $\mu$ into $\nu$ can be used to transform tests of uniform\-ity into tests of goodness of fit to $\mu$.

\subsection{Cartan--Hadamard manifolds}

Let $m$ be a point in a Cartan--Hadamard manifold, $\mathcal{X}$, and let $\mu$ and $\nu$ be probability distributions on $\mathcal{X}$ 
and $T \mathcal{X}_m$, respectively, such that the density of $\mu$ with respect to $\nu$ is positive.
By Proposition 3, there is a canonical almost-diffeo\-morphism $\phi : \mathcal{X} \to T \mathcal{X}_m$
that transforms $\mu$ into $\nu$.
Since $T \mathcal{X}_m$  can be identified with $\mathbb{R}^p$ (where $p$ is the dimension of $\mathcal{X}$),
standard goodness-of-fit tests on $\mathbb{R}^p$ can be adapted to give goodness-of-fit tests on $\mathcal{X}$.

\subsection{Simplices}

On the simplex $\Delta_{p-1}$ the uniform distribution provides a canonical choice for $\nu$. Then the transformation, $\phi$, of Proposition 4 or  Proposition 5
that takes $\mu$ into $\nu$ can be used to transform tests of uniformity into tests of goodness of fit to $\mu$.

An appealing test of uniformity on $\Delta_{p-1}$ is the score test of uniformity ($\alpha_1 = \dots = \alpha_p = 1$) within the Dirichlet family 
with densities (with respect to the uniform distribution) 
\[
 f(y_1, \ldots , y_p ; \balpha )  = \frac{\Gamma ( \sum_{j=1}^p  \alpha_j)}{\prod_{j=1}^p \Gamma (\alpha_j)} 
\prod_{j=1}^p y_j ^{\alpha_j - 1} ,
 \label{Dir}
 \]
where $\balpha = (\alpha_1, \ldots , \alpha_p)$ with $\alpha_i > 0$ for $i \in \left\{ 1, \ldots , p \right\}$.
For  independent observations $\by_1, \ldots , \by_n$ on $\Delta_{p-1}$ with 
$\by_i = (y_{i1}, \ldots , y_{ip})$ (for $i = 1, \ldots , n$),
this score test rejects uniformity for large values of 
\[
S_n =\frac{n}{\psi '(1)} 
 \left\{ \frac{\psi '(p)}{p \psi '(p) - \psi '(1)} \sum_{j=1}^p 
\sum_{k=1}^p w_j w_k
 -   \sum_{j=1}^p w_j^2  \right\} ,
\]
where $w_j = n^{-1} \sum_{i=1}^n \ln y_{ij}$ and $\psi$ denotes the digamma function.
Under uniformity the large-sample asymptotic distribution of $S_n$ is $\chi^2_p$.

\section{Simulation studies}  

\begin{figure}[h!]
\centering
\includegraphics[width = 0.85 \textwidth, height = 7 cm]{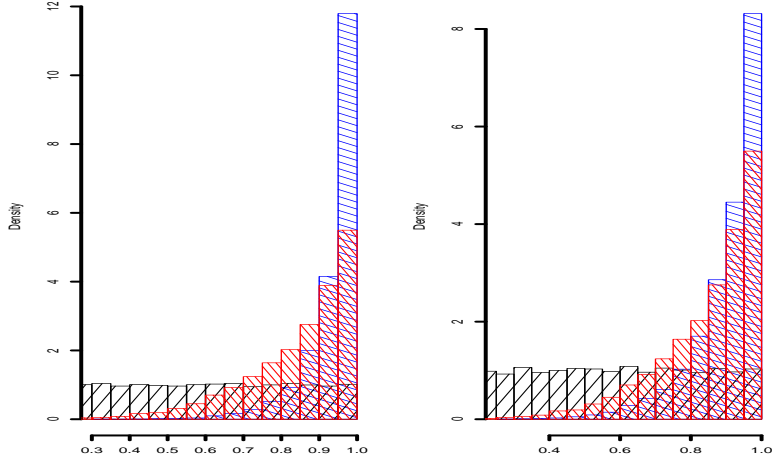}
\caption{
Behaviour of test of goodness of fit 
(a) to true $F(\bmu, 10)$ distribution on $S^2$ (black),
(b) to fitted $F(\hat{\bmu}, \hat{\kappa})$ distribution (red), 
(c) to projected normal $P \mathcal{N}_3( \bmu, \bI_3)$ distribution (blue),
using test based on Rayleigh's test of uniformity. 
The histograms are of $p$-values (based on the large-sample asymptotic $\chi^2_3$ distribution).
Each histogram summarizes 10,000 simulations, each of size $50$ (left) or $500$ (right).}
\label{fig: Sphere}
\end{figure}

In order to assess the performance of our tests, we consider  three simulation studies.
The first involves the goodness-of-fit test on $S^2$ based on the Rayleigh test of uniformity.
First $10,000$ random samples of size $50$ 
were simulated from the Fisher distribution $F(\bmu, \kappa)$ with given mode $\bmu$ and concentration $\kappa = 10$.
For each sample, goodness of fit to (a) the true $F(\bmu, 10)$ distribution,
(b) the fitted $F(\hat{\bmu}, \hat{\kappa})$ distribution, where $\hat{\bmu}$ and $\hat{\kappa}$ are the maximum likelihood 
estimates of $\bmu$ and $\kappa$, was assessed.
Then $10,000$ random samples of size $50$ 
were simulated from the projected normal $P \mathcal{N}_3( \bmu, \bI_3)$ distribution 
(obtained by projecting the trivariate normal $\mathcal{N}_3( \bmu, \bI_3)$ distribution radially onto $S^2$)
and goodness of fit to the $F(\bmu, 10)$ distribution was assessed.
The resulting $p$-values (based on the large-sample asymptotic $\chi^2_3$ distribution) 
are shown in the histograms on the left of Figure \ref{fig: Sphere}.
Corresponding histograms for $1,000$ samples of size $500$ are given on the right of Figure \ref{fig: Sphere}.
The fairly uniform distribution of $p$-values for fit to the true distribution indicates that the test detects good fit when it is present, 
whereas the clustering of $p$-values near $1$ 
when assessing goodness of fit to the fitted distribution shows the anticipated excellent fit in this case.
For samples generated from $P \mathcal{N}_3( \bmu, \bI_3)$,  the $p$-values for fit to the $F(\bmu, 10)$ distribution also cluster near $1$,
meaning that this test does not detect that the data come from the wrong model.

\begin{figure}[h!]
\centering
\includegraphics[width = 0.85 \textwidth, height = 7 cm]{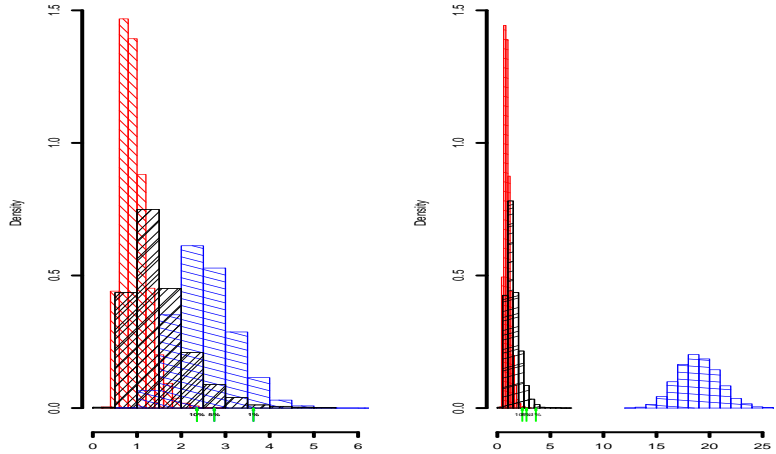}
\caption{
Behaviour of test of goodness of fit 
(a) to true $F(\bmu, 10)$ distribution on $S^2$ (black),
(b) to fitted $F(\hat{\bmu}, \hat{\kappa})$ distribution (red), 
(c) to projected normal $P \mathcal{N}_3( \bmu, \bI_3)$ distribution (blue),
using test based on Gin\'e's $F_n$ test of uniformity. 
The histograms are of values of $F_n$.  
Each histogram summarizes 10,000 simulations, each of size $50$ (left) or $500$ (right). 
Green arrows on horizontal axes are 10\%, 5\% and 1\% upper quantiles of asymptotic distribution.}   
\label{fig: Gine}
\end{figure}

One possible explanation for the inability of the above test to detect that the data come from the wrong model is that the Rayleigh test 
of uniformity is not consistent against all alternatives.
Therefore a second simulation study was carried out, which was like the first but with the Rayleigh test replaced by Gin\'e's \cite{Gine75} $F_n$ test
\cite[Section 10.4.1]{MardiaJupp00}, 
which is consistent against all alternatives to uniformity on $S^2$.
Histograms of the resulting values of $F_n$ are shown in Figure \ref{fig: Gine} for sample sizes, $n$, of $50$ (left) and $500$ (right).
Significance was assessed using the asymptotic quantiles given in \cite{Keilson.etal83}  \cite[Section 10.4.1]{MardiaJupp00}.
For assessing goodness of fit to the true distribution,
the proportions of the values of the statistic that exceeded the asymptotic 10\%, 5\% and 1\% upper quantiles were 
$0.10, 0.05$ and $0.01$ ($n= 50$) and
$0.10, 0.04$ and $0.01$ ($n= 500$), respectively,
indicating that the test detects good fit when it is present. 
For fit to the fitted distribution, none
of the values of $F_n$ exceeded the asymptotic 10\% quantile, indicating the anticipated excellent fit in this case.
For samples generated from $P \mathcal{N}_3( \bmu, \bI_3)$,
the proportions of the values of $F_n$ that exceeded the asymptotic 10\%, 5\% and 1\% upper quantiles were 
$0.58, 0.34$ and $0.05$ for $n= 50$,
while for $n= 500$, all the values of $F_n$ far exceeded the asymptotic 1\% upper quantile.
This indicates clearly that the the test can detect bad fit. 

\begin{figure}[h!]
\centering
\includegraphics[width = 0.85 \textwidth, height = 7 cm]{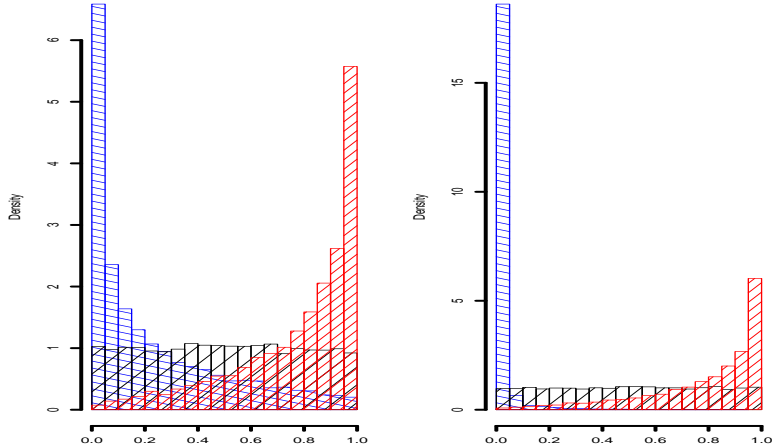}
\caption{
Behaviour of test of goodness of fit
(a)  to true isotropic Mardia--Dryden $MD( [ \bmu ], 0.125)$ distribution on $\Sigma_2^5$  (black),
(b)  to fitted isotropic $MD( [ \hat{\bmu} ], \hat{\kappa})$ distribution (red), 
(c) to non-isotropic Mardia--Dryden distribution obtained by Gaussian 
$\mathcal{N}_2( \bzero , \mathrm{diag} (1, 25))$  perturbations of $\bmu$ (blue),
using test based on Mardia's test of uniformity.
The histograms are of $p$-values (based on the large-sample asymptotic $\chi^2_{15}$ distribution).
Each histogram summarizes 10,000 simulations, each of size $50$ (left) or $500$ (right).}   
\label{fig: shape}
\end{figure}

The third simulation study involves the goodness-of-fit test on $\Sigma^5_{2}$ based on 
Mardia's  \cite{Mardia99} test of uniformity. 
First, $10,000$ random samples of size $50$ 
were simulated from the isotropic Mardia--Dryden $MD( [ \bmu ], 0.125)$ distribution with given mode $[ \bmu ]$.
For each sample, goodness of fit to    \linebreak
(a) the true $MD( [ \bmu ], 0.125)$ distribution,
(b) the fitted $MD( [ \hat{\bmu} ], {\hat \kappa})$ distribution, where $[ \hat{\bmu} ]$ and $\hat{\kappa}$ are the maximum likelihood 
estimates of $[ \bmu ]$ and $\kappa$
(calculated by the EM method of \cite{KumeWelling10}), 
was assessed using Mardia's uniformity test on $\Sigma^5_2$.
Then $10,000$ random samples of size $50$ 
were simulated from the non-isotropic Mardia--Dryden distribution obtained by Gaussian 
$\mathcal{N}_2( \bzero , \Sigma)$  perturbations of $\bmu$,
where $\Sigma = \mathrm{diag} (1, 25)$, and goodness of fit to the 
$MD( [ \bmu ], 0.125)$ distribution was assessed.
The resulting $p$-values 
based on the large-sample asymptotic $\chi^2_{15}$ distribution 
are shown in the histograms on the left of Figure \ref{fig: shape}.
Corresponding histograms for $10,000$ samples of size $500$ are given on the right.
The fairly uniform distribution of $p$-values for fit to the true distribution indicates that the test detects good fit when it is present. 
The clustering of $p$-values near $1$ for fit to the fitted distribution shows the anticipated excellent fit in this case.
For samples generated from the non-isotropic distribution, the $p$-values cluster near $0$, 
indicating that the test can detect bad fit.

\begin{thebibliography}{99}
\bibitem{Beran79}
Beran, R.  (1979).
Exponential models for directional data. 
\emph{Ann. Statist.} {\bf 7},  1162--1178.
\bibitem{BestFisher86}
 Best, D.J. \& Fisher,  N.I.  (1986).
Goodness-of-fit and discordancy tests for samples
from the Watson distribution on the sphere.
\emph{ Austral. J. Statist.} {\bf 28},  13--"31.
\bibitem{BoulDuch97}
Boulerice, B.  \& Ducharme,  G.  (1997).
Smooth tests of goodness-of-fit for directional and axial data.
\emph{J. Multivariate Anal.} {\bf 60}, 154--175.
\bibitem{DrydenMardia16}
Dryden, I.L.  \& Mardia, K.V. (2016).
\emph{Statistical Shape Analysis with applications in R}, second ed.
Wiley, Chichester.
\bibitem{FisherBest84}
Fisher,  N.I. \& Best,  D. (1984).
Goodness-of-fit tests for Fisher's   \linebreak
distribution on the sphere. %
\emph{Austral. J. Statist.} {\bf 26}, 142--150.
\bibitem{Gallotetal93}
Gallot, S., Hulin, D. \&  Lafontaine, J. (1993).
\emph{Riemannian Geometry},
second ed. 
Springer-Verlag, Berlin.
\bibitem{Gine75}
Gin\'e, R. (1975).
Invariant tests for uniformity on compact Riemann\-ian 
manifolds based on Sobolev norms.
\emph{Ann. Statist.} {\bf 3}, 1243--1266.
\bibitem{Helgason78}
Helgason, S.  (1978).
\emph{Differential Geometry, Lie Groups and  Symmetric Spaces.} 
Academic Press, New York.	     
\bibitem{Huckemannetal10}
Huckemann,  S.F.,  Kim,  P.T., Koo, J.Y.  \& Munk, A.  (2010).
M\"{o}bius deconvolution on the hyperbolic plane with application to impedance density 
 estimation.
\emph{Ann. Statist.} {\bf 38}, 2465--2498.
\bibitem{Jonesetal15}
Jones, M.C.,  Pewsey, A.  \&  Kato,  S. (2015).
On a class of circulas for circular distributions.
\emph{Ann. Inst. Statist. Math.} {\bf 67}, 843--862.
\bibitem{Jung.etal12}
Jung, S., Dryden, I.L. \& Marron, J.S. (2012).
Analysis of principal nested spheres. 
\emph{Biometrika} {\bf 99}, 551--568.
\bibitem{Jupp05}
Jupp, P.E.  (2005).
Sobolev tests of goodness of fit of distributions on compact Riemannian manifolds.  
\emph{Ann. Statist.} {\bf 33},  2957--2966.
\bibitem{Jupp15}
Jupp, P.E.  (2015).
Copulae on products of compact Riemannnian \linebreak
manifolds.   
\emph{J. Multivariate Anal.} {\bf 14}, 92--98.
\bibitem{JuppSpurr83}
Jupp, P.E.   \& Spurr, B.D. (1983).
Sobolev tests for symmetry of   \linebreak
directional data.   
\emph{Ann. Statist.} {\bf 11}, 1225--1231.
\bibitem{JuppSpurr85}
Jupp, P.E.   \& Spurr, B.D. (1985).
Sobolev tests for independence of directions. 
\emph{Ann. Statist.} {\bf 13}, 1140--1155.
\bibitem{Keilson.etal83}
Keilson, J., Petrondas, D., Sumita, U. \& Wellner, J. (1983).
Significance points for some tests of uniformity on the sphere.
\emph{J. Statist. Comput. Simulation} {\bf 17}, 195--218.
\bibitem{Kendalletal99}
Kendall, D.G., Barden, D., Carne, T.K. \& Le, H. (1999).
\emph{Shape and Shape Theory.}
Wiley, Chichester.
\bibitem{Kent82}
Kent,  J.T. (1982).
The Fisher--Bingham distribution on the sphere.
\emph{J. R. Stat. Soc.} {\bf B  44},  71--180. 
\bibitem{KobayashiNomizu69}
Kobayashi, S. \& Nomizu,  K. (1969).
\emph{Foundations of Differential      \linebreak 
Geometry,    Volume II}. 
Interscience, New York.
\bibitem{KumeWelling10}
Kume, A. \& Welling, M. (2010).
Maximum likelihood estimation for the offset-normal shape distributions using EM.
\emph{J. Comp. Graphical Stat.} {\bf 19}, 702--723.
\bibitem{Le.Small99}
Le, H.   \& Small, C.G. (1999).
Multidimensional scaling of simplex shapes.   
\emph{Pattern Recognition} {\bf 32}, 1601--1613.
\bibitem{Mardia99}
Mardia, K.V. (1999).
Directional statistics and shape analysis, 
\emph{J. Appl. Statist.} {\bf  26}, 949--957.
\bibitem{MardiaJupp00}
Mardia, K.V.  \& Jupp, P.E. (2000).
\emph{Directional Statistics}.
Wiley, \linebreak
Chichester.   
\bibitem{Mardiaetal84}
Mardia, K.V., Holmes, D. \& Kent, J.T. (1984).
A goodness-of-fit test for the von Mises--Fisher distribution.
\emph{J. R. Stat. Soc.} {\bf B 46},  72--78.
\bibitem{Moser65}
Moser, J. (1965).
On the volume elements on a manifold.
\emph{Trans. Amer. Math. Soc.} {\bf 120}, 286--294.
\bibitem{PewseyKato16}
Pewsey, A. \& Kato, S. (2016).
Parametric bootstrap goodness-of-fit testing for Wehrly--Johnson bivariate circular distributions,
\emph{Statistics and Computing} 
{\bf 26}, 1307--1317.    
\bibitem{Rosenblatt52}
Rosenblatt, M. (1952).
Remarks on a multivariate transformation.
\emph{Ann. Math. Statist.} {\bf 23}, 470--1472.
\bibitem{Small96}
Small, C.G. (1996).
\emph{The Statistical Theory of Shape}.
Springer, New York.
\end{thebibliography}
\end{document}